\documentclass[a4paper, 11pt]{article}

\usepackage[english]{babel}
\usepackage[all]{xy}
\usepackage{graphicx}
\usepackage{amsmath,amsthm}
\usepackage{amsfonts,amssymb}
\usepackage{float}
\usepackage{bbm}

\newcommand{\rz}[1]{\mathbb{#1}}
\def\N{\rz N}  
\def\E{\rz E}  
\def\Q{\rz Q}  

\def\Pb{\rz P}
\def\shf{{\cal F}}
\def\she{{\cal E}}

\newcommand{\indi}{\mathbbm{1}}

\newtheorem{theo}{Theorem}[section]
\newtheorem{lem}[theo]{Lemma}
\newtheorem{prop}[theo]{Proposition}

\newtheorem{defi}[theo]{Definition}

\newtheorem{remark}[theo]{Remark}

\begin{document}

\title{\textbf{\Large Mathematical model for resistance and optimal strategy}}
 \author{\small{
 Blandine B\'ERARD BERGERY\footnote{IECN, Nancy-Universit\'e, CNRS, INRIA, Boulevard des Aiguillettes B.P. 239 F-54506 Vand\oe uvre l\`es Nancy},
 Christophe PROFETA$^*$,
 Etienne TANR\'E\footnote{INRIA, EPI TOSCA, 2004 route des Lucioles, BP93, 06902 Sophia-Antipolis, France}}}
\date{\today}
\maketitle

\hrule

\smallskip
\textbf{Abstract:} We propose a mathematical model for one pattern of charts studied in technical analysis: in a phase of consolidation, the price of a risky asset goes down $\xi$ times after hitting a resistance level.  We construct a mathematical strategy and we calculate the expectation of the wealth for the logaritmic utility function. Via simulations, we compare the strategy with the standard one.


\vskip 3mm

\hrule

\section{Introduction}\label{sec0}

\textbf{1.} In financial market, investors can use three approaches to decide on their future investment: the fundamental approach, the mathematical approach and the technical analysis. The last one is popular among traders and is based upon the study of the charts and the past behavior of the price  .

The main hypothesis of technical analysis is that all the needed information is contained inside the records of prices and transaction volumes. The analysis of the charts is therefore sufficient to predict the near future of the price. This hypothesis contradicts most of mathematical models, which are essentially markov. Then, technical analysis seems to have no theoretical justifications, and even no proven efficiency (see \cite{a1}). Recent works as \cite{a3,i2,martinezetal} tries to compare this approach to the mathematical one.

Let us consider a classic pattern of charts studied via technical analysis. In a phase of consolidation, a price does not have either ascending trend, nor decreasing trend. The price moves between two barrier levels: the upper one is called resistance, the lower one is called support (see \cite{b6} for precise statement). When the price bounces three times on the support, price will likely not go down the bottom barrier. We can expect a rise and that the price will go up through the upper barrier. This belief will have an influence on trader's behavior, and therefore on price. However, this kind of rule has no mathematical justification.

The standard Black and Scholes model does not show this kind of behavior, although observations of real charts show phenomenon of support and resistance. In this article, we propose a mathematical model derived from the Black and Scholes model which has  resistance. The trajectory makes a random number of downcrossings between two fixed levels before  leaving. Then, we evaluate the optimum portfolio strategy in the case of a logarithmic utility function. The efficiency of the strategy is compared with the classic one.

\textbf{2.} First, we recall the classic problem of portfolio optimization. We consider a trader whose portofolio is composed of 2 assets:
\begin{itemize}
  \item a risky asset $S$ solution of the stochastic differential equation (SDE):
\begin{equation}\label{S}  
\frac{dS_t}{S_t}=\mu_0 dt+\sigma dB_t,
\end{equation}
  \item and a bond $S^0$ whose price evolves according to: $\frac{dS^0_t}{S^0_t}=r dt.$
 \end{itemize}
Then, the trader's wealth at time $t$ is given by:
\begin{equation}
\label{Wbis} 
\frac{dW^\pi_t}{W^\pi_t}=\pi_t \sigma dB_t+ \left(r+\pi_t (\mu_0-r)\right)dt.
\end{equation}
where $\pi_t$ stands for the proportion of risky asset $S$ the trader holds at time $t$. The aim of the trader is to determine the strategy $\pi^c$ which maximizes $\E[\log(W^\pi_T)]$, where $T$ is a fixed time. Then, the optimal strategy $\pi^c$ is:
\begin{equation}
\label{strategieinitial}  
\pi^c_s=\frac{\mu_0-r}{\sigma^2}.
\end{equation}

\begin{remark}
In this initial problem, the risky asset $S$ is a geometric brownian motion, which can be written under the form:
\begin{equation}
\label{solutions} 
S_t=S_0\exp\left(\sigma B_t + (\mu_0-\frac{\sigma^2}{2})t\right)=S_0\exp\left(\sigma \left(B_t+\mu t\right)\right)
= S_0\exp\left(\sigma B_t^\mu \right),
\end{equation}
where $\displaystyle\mu:=\frac{1}{\sigma}\left(\mu_0-\frac{\sigma^2}{2}\right)$ and $B_t^\mu :=B_t+\mu t$. As a result, studying only the moves of $B^\mu$ is enough to get the moves of $S$. We will almost only consider $B^\mu$ in the article, with $\mu \neq 0$.
\end{remark}

\textbf{3.}  The risky asset $S$ in item \textbf{2} does not present support or resistance levels. The aim of the following sections is to propose a price process derivated from $S$ which has resistance level. In charts, a resistance levels is drawn thanks to aligned local maxima observed on the past trajectory. The observations are obviously not perfectly exact, thus we consider that the resistance line is thick. We fix three levels: $S_0^-$, $S_0$ and $S_0^+$, where $[ S_0 ; S_0^+]$ depicts the resistance line drawn in charts.  We first construct a process $X$ derivated from $B^\mu$. Then, the price process is $Z=S_0e^{\sigma X}$ and $Z$  does at least $n$ downcrossings from level $S_0$ to a fixed level $S_0^-$ before being allowed to reach the upper bound of the resistance $S_0^+$. The model would be useless if we were not able to compute a strategy. With our model, we get an explicit formula for the optimal strategy and thus we make numerical simulations.

The study has four steps:
\begin{itemize}
  \item study of the downcrossings of  $B^\mu$ in Section \ref{sec1}, which is the main tool for all our results,
    \item construction of  $X$ and $Z$  in Section \ref{sec1bis},
  \item computation of the optimal strategy if the number of downcrossings  is a fixed integer $n$ in Section \ref{sec2},
  \item computation of the optimal strategy if the number of downcrossings  is a random variable $ \xi$ in Section \ref{sec3},
  \item comparison of our model with the classic one throughout numerical simulations in Section \ref{sec4}.
\end{itemize}


\section{Study of the downcrossings of a brownian motion with drift} \label{sec1}

We consider a Brownian motion with drift $\mu\neq 0$  $(B_t^\mu:=B_t+\mu t, t\geq0)$ started from 0, $(\shf_t)_{t\geq0}$ its natural filtration, and  two fixed reals $\epsilon>0$ and $\alpha>0$.  First, we define and study the event $A_n$: ``$B^\mu $ has done at least $n$ downcrossings from 0 to $-\alpha$ before hitting the level $\epsilon$" in item \textbf{1}. Then we study in item \textbf{2} the martingale $M_t^n:=\Pb(A_n|\shf_t)$. Finally, we study $B_t^\mu$ conditioned by $A_n$ in item \textbf{3}.

\textbf{1.}
To begin with, let us consider the sequence of stopping times at levels 0 and $-\alpha$:
\[
\begin{array}{rclrcl}
\sigma_1&=&\inf\{t\geq0; B_t^\mu=-\alpha\},&
\sigma_2&=&\inf\{t\geq\sigma_1; B_t^\mu=0\}\\
\sigma_{2k+1}&=&\inf\{t\geq\sigma_{2k}; B_t^\mu=-\alpha\},&
\sigma_{2k+2}&=&\inf\{t\geq\sigma_{2k+1}; B_t^\mu=0\}.
\end{array}
\]
The $k+1^\text{th}$ downcrossing from $0$ to $-\alpha$ takes place between the times $\sigma_{2k}$ and $\sigma_{2k+1}$. $\sigma_{2k+2}$ is the time where $B$ comes back to 0 after the $k+1^\text{th}$ downcrossing. For example, $k=0$ corresponds to the first downcrossing. %
Let us remark that, if the process does not reach $-\alpha$, we consider nevertheless that it is in a downcrossing phase, from time 0 to infinite time. To take this event into account, we use the convention $  \sigma_{-1} = 0$.

Schematically, we have roughly this pattern:
$$  
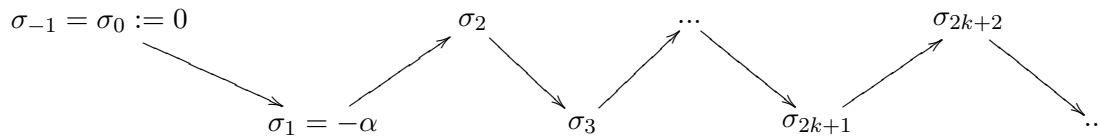
\xymatrix{
  \sigma_{-1}=\sigma_0:=0 \ar[rd]& & \sigma_2 \ar[rd]& & ... \ar[rd] & & \sigma_{2k+2} \ar[rd] \\
     & \sigma_1= -\alpha \ar[ru] & & \sigma_3 \ar[ru] & & \sigma_{2k+1} \ar[ru] & & ...
} $$
Finally, we define $T_\epsilon:=\inf\{t\geq0, B_t^\mu=\epsilon\}$ the hitting time of $\epsilon$ by $B^\mu$, which is the ending time of our study.

The event $A_n$: ``the trajectory has done at least $n$ downcrossings from 0 to $-\alpha$ before hitting the level $\epsilon$" can be expressed through these stopping times:
\begin{equation}
\label{an} 
A_n:=\{\sigma_{2n-1}< T_\epsilon\}
\end{equation}
for $n\geq0$. By definition, we have the inclusions: $ A_0 \supset A_1 \supset \dots \supset A_n \supset \dots$.

Then, the risky asset will be the process $(X_t:=B_t^\mu, 0\leq t )$ conditioned by $A_n$, $n\geq1$. The aim of the two following items is to establish the equations satisfied by $X$.

Let us start our study of $A_n$ by computing $\Pb(A_n)$.

\begin{lem}\label{probaAn}
$$\Pb(A_n)= p^n \quad \textrm{ where } p=\frac{S(\epsilon)-S(0)}{S(\epsilon)-S(-\alpha)}  \Pb_{-\alpha}( T_0 < \infty),$$
and $S(x):=e^{-2\mu x}$ is a scale function for $B^\mu$. If $\mu \geqslant 0$, then  $ \Pb_{-\alpha}( T_0 < \infty) = 1$, else $ \Pb_{-\alpha}( T_0 < \infty) = e^{2\mu \alpha}$.
\end{lem}

\textbf{Proof.} If $n=0$, we have directly $\Pb(A_0)=\Pb(\sigma_{-1}<T_\epsilon)=1$.

For $n\geq1$, we use the successive inclusions of the $A_i$ to write:
$$\{\sigma_{2n-1}<T_\epsilon\} = \{\sigma_1<T_\epsilon\}\cap\{\sigma_3<T_\epsilon\}\cap...\cap\{\sigma_{2n-1}<T_\epsilon\}.$$
Then, we get by conditioning by $A_1$:
\[
\Pb (A_n) = \E(  \E( \indi_{\{\sigma_{2n-1}<T_\epsilon\}}| A_1) )= \E( \indi_{\{\sigma_1<T_\epsilon\}} \E( \indi_{\{\sigma_3<T_\epsilon\}\cap...\cap\{\sigma_{2n-1}<T_\epsilon\}}| A_1) )  .
\]
Under $A_1$, we have to consider two cases: $\sigma_2 < \infty$ and $\sigma_2 = \infty$.
\begin{itemize}
\item  If $\sigma_2 = \infty$, then $\sigma_k, k > 2$ and $T_\epsilon$ are infinite and $\E(  \indi_{\{\sigma_3<T_\epsilon\}\cap...\cap\{\sigma_{2n-1}<T_\epsilon\}}| A_1 \cap \{\sigma_2 = \infty\}) ) =0$.
\item If $\sigma_2 < \infty$,  since $\sigma_1 < T_\epsilon$, $B_\mu$ must cross the level $0$ to go up to the level $\epsilon$. This means that $\sigma_2 \leqslant T_\epsilon$. Since $B_{\sigma_2}^\mu =0$, the process $(\overline{B}^\mu_t =  B^\mu_{t+\sigma_2}  - B_{\sigma_2}^\mu$, $t\geqslant 0)$,  is a brownian motion with drift $\mu$  independent from $(B_{t}^\mu, t\in [0, \sigma_2])$. $B^\mu$ having already done one downcrossing, we have to remove one downcrossing from those of $\overline{B}^\mu$:
$$\E( \indi_{\{\sigma_3<T_\epsilon\}\cap...\cap\{\sigma_{2n-1}<T_\epsilon\}}| A_1 \cap \{\sigma_2 < \infty\}) = \E( \indi_{\{\overline{\sigma}_1<T_\epsilon\}\cap...\cap\{\overline{\sigma}_{2n-3}<T_\epsilon\}}).$$
\end{itemize}
Let us note that $\Pb( \sigma_2 < \infty |A_1) = \Pb_{-\alpha}( T_0 < \infty)$, probability to hit 0 starting from $-\alpha$. Hence
$$\Pb (A_n) = \E( \indi_{\{\sigma_1<T_\epsilon\}} )  \Pb_{-\alpha}( T_0 < \infty)  \E( \indi_{\{\overline{\sigma}_1<T_\epsilon\}\cap...\cap\{\overline{\sigma}_{2n-3}<T_\epsilon\}}).
$$
Repeating this method, we have by induction $\Pb (A_n) = \left[ \Pb(\sigma_1<T_\epsilon)  \Pb_{-\alpha}( T_0 < \infty) \right]^n$. The quantity $\Pb(\sigma_1<T_\epsilon)$ is the probability to hit the level $-\alpha$ before the level $\epsilon$, starting from 0. Consequently,
\[
\Pb (A_n) =\left(\frac{S(\epsilon)-S(0)}{S(\epsilon)-S(-\alpha)}  \Pb_{-\alpha}( T_0 < \infty)\right)^n,
\]
where $S(x)$ a scale function for $B^\mu$. $ \Pb_{-\alpha}( T_0 < \infty) = 1$ if $\mu \geqslant 0$ and $ \Pb_{-\alpha}( T_0 < \infty) = e^{2\mu \alpha}$ if $\mu < 0$ (see \cite{b1}).
\qed

\textbf{2.} Next, we study the martingale $M_t^n:=\Pb(A_n|\shf_t)$ in the 2 following lemmas. This martingale will play a central role when conditioning by $A_n$ (c.f. Proposition \ref{Btilde}).
\begin{lem}\label{defM}
Let us consider the martingale \(M_t^n:=\Pb(A_n|\shf_t)\).
\begin{multline*}
M_t^n
=\indi_{\{\sigma_{2n-1}<t\wedge T_\epsilon\}}
+\indi_{\{t<T_\epsilon\}}\sum\limits_{i=0}^{n-1} \indi_{\{\sigma_{2i}\leq t < \sigma_{2i+1}\}} p^{n-1-i} \frac{S(\epsilon)-S(B_t^\mu)}{S(\epsilon)-S(-\alpha)}\quad\text{(Downcrossings)}\\
 + \indi_{\{t<T_\epsilon\}}\sum\limits_{i=1}^{n-1} \indi_{\{\sigma_{2i-1}\leq t <\sigma_{2i}\}} p^{n-i}.\quad\text{(Upcrossings)}
\end{multline*}
\end{lem}
\textbf{Proof.} We start by dividing the event $A_n$ in a disjointed union of events :
$$A_n=\left(\bigcup\limits_{k=0}^{2n-2} (A_n\cap\{\sigma_k \leq t <\sigma_{k+1}\})\right) \cup (A_n\cap\{t\geq\sigma_{2n-1}\}).$$
Let us remark first that the last term can be rewritten: $A_n\cap\{t\geq\sigma_{2n-1}\}=\{\sigma_{2n-1}<t\wedge T_\epsilon\}$ and is thus $\shf_t$-mesurable. Consequently:
\begin{equation}\label{M}
M_t^n=\indi_{\{\sigma_{2n-1}<t\wedge T_\epsilon\}}+\sum\limits_{k=0}^{2n-2} \Pb(A_n\cap\{\sigma_k \leq t <\sigma_{k+1}\}|\shf_t).
\end{equation}
Let us define $C_k=A_n\cap\{\sigma_k \leq t <\sigma_{k+1}\}$. We must now separate the cases $k$ is even and $k$ is odd.

\begin{itemize}
\item[i)] For $k=2i$, $0\leq i \leq n-1$:
$$C_{2i}=\{\sigma_{2i}\leq t <\sigma_{2i+1}\}\cap\{\sigma_{2n-1}< T_\epsilon\}.$$
At time $t$, the process has still $n-i$ downcrossings to the level $-\alpha$ to do before 
hitting 
the level $\epsilon$.
 The first downcrossing takes place from $B_t^\mu$ to $-\alpha$, and the $(n-i-1)^\text{th}$ following ones from 0 to $-\alpha$. 
Computations similar to Lemma~\ref{probaAn} gives
\begin{equation}\label{descentes}
\Pb(C_{2i}|\shf_t)=\indi_{\{\sigma_{2i}\leq t<\sigma_{2i+1},\; t<T_\epsilon \}} \frac{S(\epsilon)-S(B_t^\mu)}{S(\epsilon)-S(-\alpha)}\left(\frac{S(\epsilon)-S(0)}{S(\epsilon)-S(-\alpha)} \Pb_{-\alpha} (T_0 < \infty) \right)^{n-i-1}.
\end{equation}
\item[ii)] For  $k=2i-1$, $1\leq i \leq n-1$,
\[
C_{2i}=\{\sigma_{2i-1}\leq t <\sigma_{2i}\}\cap\{\sigma_{2n-1}\leq T_\epsilon\}.
\]
At time $t$, the process has still $n-i$ downcrossings to the level $-\alpha$ to do before 
hitting the level $\epsilon$. These downcrossings will all take place from 0 to $-\alpha$, since the process must first reach the level 0 (at $\sigma_{2i}$) before the next downcrossing to be taken into account. Then we get:
\begin{equation}\label{montees}
\Pb(C_{2i-1}|\shf_t)=\indi_{\{\sigma_{2i-1}\leq t<\sigma_{2i},\; t<T_\epsilon \}} \left(\frac{S(\epsilon)-S(0)}{S(\epsilon)-S(-\alpha)} \Pb_{-\alpha} (T_0 < \infty) \right)^{n-i}.
\end{equation}
\end{itemize}
Finally, injecting (\ref{descentes}) and (\ref{montees}) in (\ref{M}) ends this proof.
\qed

\begin{lem}\label{dM}
The function $t\mapsto M_t^n$ is continuous and
$$dM_t^n = \sum\limits_{i=0}^{n-1} \indi_{\{\sigma_{2i}\leq t < \sigma_{2i+1}\}} 2\mu p^{n-1-i} \frac{S(B_t^\mu)}{S(\epsilon)-S(-\alpha)}\indi_{\{t< T_\epsilon\}} dB_t.$$
\end{lem}
\textbf{Proof.} We use the expression of $M^n$ given by Lemma \ref{defM} on each interval $[\sigma_k, \sigma_{k+1}[$, $-1 \leqslant k \leqslant 2n$. Since $S(B_t^\mu)=\exp\left(-2\mu(B_t+\mu t)\right)$, applying It\^{o}'s formula gives $dS(B_t^\mu)=-2\mu S(B_t^\mu)dB_t$ and the result holds.
\qed

\textbf{3.} We end this section with the main result, concerning the decomposition of the Brownian motion $B$ under  the event $A_n$, which will enable us to condition by $A_n$.
\begin{prop}\label{Btilde}
There exists a Brownian motion $\widetilde{B}$ independent from $A_n$ such that, on the event $A_n$:
$$B_t=\widetilde{B}_t+2\mu \int_0^t \left(\sum\limits_{i=0}^{n-1} \indi_{[\sigma_{2i},\sigma_{2i+1}[}(s) \frac{S(B_s^\mu)}{S(\epsilon)-S(B_s^\mu)}\right)ds.$$
\end{prop}
\textbf{Proof.} We use formula (1)  from \cite{i1} (Chapter II p.45).  Let $\shf_t^n:=\shf_t\vee \sigma(A_n)$. On the event $A_n$, we have:
\begin{equation}\label{Yor}
B_t=\widetilde{B}_t+\int_0^t \frac{d<B,M^n>_s}{M_s^n},
\end{equation}
where $\widetilde{B}$ is a $\shf_t^n$ Brownian motion. Since $\widetilde{B}_t$ is independent from $\shf_0^n$, it is independent from $A_n$. Then, we deduce from Lemma \ref{dM}:
\[
d\left<B,M^n\right>_s=\sum\limits_{i=0}^{n-1} \indi_{\{\sigma_{2i}\leq s < \sigma_{2i+1}\}} 2\mu p^{n-1-i} \frac{S(B_s^\mu)}{S(\epsilon)-S(-\alpha)}\indi_{\{s< T_\epsilon\}} du.
\]
From Lemma \ref{defM}, $\displaystyle M_s^n=p^{n-1-i} \frac{S(\epsilon)-S(B_s^\mu)}{S(\epsilon)-S(-\alpha)}$ for $\{\sigma_{2i}\leq s < \sigma_{2i+1}\}$. As a result,  (\ref{Yor}) becomes:
$$B_t=\widetilde{B}_t+2\mu \int_0^t \left(\sum\limits_{i=0}^{n-1} \indi_{\{\sigma_{2i}\leq s<\sigma_{2i+1}\}} \frac{S(B_s^\mu)}{S(\epsilon)-S(B_s^\mu)}\right)\indi_{\{s< T_\epsilon\}} ds.$$
But, on $A_n$,  $s\leq \sigma_{2i+1}$ implies $ s< T_\epsilon$, which gives the result.
\qed

\section{Construction of the price process}\label{sec1bis}

 Let $S_0$ be the starting price of the risky asset and $S_0^-$ , $S_0^+$ be two fixed levels with $0<S_0^-<S_0<S_0^+$. The Section describes the construction of a price process $Z$ such that $Z$ does $n$ downcrossings from level $S_0$ to level $S_0^-$ before being allowed to reach the upper bound of the resistance line $S_0^+$. Let us note that the number of downcrossings $n$ may be random. The definition of $Z$ is given in item  \textbf{1}, then $Z$ is linked to the classic price $S$ in item \textbf{2}. Finally, in item \textbf{3}, we briefly propose another construction  via local time.

\textbf{1.} To begin with, we define the underlying process $X$. We denote  $\displaystyle\alpha=-\frac{1}{\sigma}\log\left(\frac{S_0^-}{S_0}\right)$ and $\displaystyle\epsilon=\frac{1}{\sigma}\log\left(\frac{S_0^+}{S_0}\right)$.

\begin{defi}
\label{definitionX}
The process $X$ is defined  for all $t\in [0,T_\epsilon]$ as follows.\\
Let
$$
\left\{
\begin{array}{rclrcl}
\sigma_1&=& \inf\{t\geq0; X_t=-\alpha\}\ ,& \sigma_2&=&\inf\{t\geq\sigma_1; X_t=0\}, \\
\sigma_{2k+1}&=&\inf\{t\geq\sigma_{2k}; X_t=-\alpha\}\ , &\sigma_{2k+2}&=&\inf\{t\geq\sigma_{2k+1}; X_t=0\},\\
T_\epsilon&=&\inf\{t\geq0, X_t=\epsilon\}.&&&
\end{array}
\right.
$$
If
\begin{itemize}
\item  $t \geqslant  \sigma_{2n -1}$ , then $X$ is a standard Brownian motion with drift $\mu$ starting at $-\alpha$ at time $\sigma_{2n-1}$,
  \item if $t\in [\sigma_{2i -1}, \sigma_{2i }]$, $i \in [0, n-1] $, then $X$ is a standard Brownian motion with drift $\mu$ starting at $-\alpha$ at time $\sigma_{2i -1}$,
   \item if $t\in [\sigma_{2i }, \sigma_{2i +1}]$, $i \in [0, n-1]$, then $X$ is solution of the equation
   \begin{equation}
\label{xdescent}
dX_t = d\widetilde{B}_t-\mu \coth(\mu(\epsilon-X_t))dt,
\end{equation}
where $\widetilde{B}$ is a standard Brownian motion.
\end{itemize}
\end{defi}

Then, the price process is $Z=S_0e^{\sigma X}$. In the next item, we explain how $Z$ (resp. $X$) is related to the standard price process $S$ (resp. the Brownian motion with drift $B^\mu$). This relation will allow us to compute explicit strategy in Sections \ref{sec2} and \ref{sec3}.

\textbf{2.} We consider the risky asset $S$ solution of the classic SDE (\ref{S}). We recall that the solution of (\ref{S}) is $S_t=S_0\exp\left(\sigma B_t^\mu \right)$, with $\mu:=\frac{1}{\sigma}\left(\mu_0-\frac{\sigma^2}{2}\right)$ and $B_t^\mu$ the Brownian motion with drift studied in Section \ref{sec1}. The wealth $W^\pi_t$ the trader holds at time $t$ is given by the equation (\ref{Wbis}) where $\pi_t$ is the proportion of the risky asset $S$ the trader holds at time $t$.

The levels $S_0^-,S_0,S_0^+$ of $S$ are the ones of $Z$ and we have the following equivalences:
\begin{enumerate}
  \item $S_t = S_0 \Leftrightarrow B_t^\mu = 0$,
  \item $S_t = S_0^- \Leftrightarrow B_t^\mu = -\alpha$ with  $\displaystyle\alpha=-\frac{1}{\sigma}\log\left(\frac{S_0^-}{S_0}\right)$,
  \item $S_t = S_0^+ \Leftrightarrow B_t^\mu = \epsilon $ with $\displaystyle\epsilon=\frac{1}{\sigma}\log\left(\frac{S_0^+}{S_0}\right)$.
\end{enumerate}
So, $T^+$ the hitting time of level $S_0^+$ by $S$ equals $T_\epsilon$, which is defined at the beginning of Section \ref{sec1}.  Similarly, the downcrossings of the process $S$ from level $S_0$ to level $S_0^-$ are exactly the downcrossings of the process $B_t^\mu$ from level 0 to level $-\alpha$. Finally, the event ``the process $S$ has done at least $n$ downcrossings from level $S_0$ to level $S_0^-$ before hitting the level $S_0^+$ " is the same as the event $A_n=\{\sigma_{2n-1}< T_\epsilon\}$, which was defined and studied in Section \ref{sec1}.

\begin{prop}
The law of the process $Z$ (resp. $X$) is  the same as the law of $S$ (resp. $B^\mu $)  conditioned by $A_n$. We denote $\Q (.)= \Pb( .| A_n)$.
\end{prop}

Thanks to this result,  the expected wealth for the logarithmic utility function is $\Pb( W^\pi | A_n)$ and can be explicitly calculated. Thus, we can evaluate the optimal strategy (see Section \ref{sec2} and \ref{sec3}).

\textbf{Proof.} Let us note $X=B^\mu$ conditioned by $A_n$. Using Proposition \ref{Btilde}, we obtain an explicit description of $X$ on $A_n$:
$$X_t= B_t + \mu t=\widetilde{B}_t+2\mu \int_0^t \left(\sum\limits_{j=0}^{n-1} \indi_{[\sigma_{2j},\sigma_{2j+1}[}(s) \frac{S(B_s^\mu)}{S(\epsilon)-S(B_s^\mu)}\right)ds + \mu t, \quad \forall t < T_\epsilon.$$
Let us consider the two cases:
\begin{itemize}
  \item  If $t \in  [\sigma_{2i -1}, \sigma_{2i }[$, then
\begin{eqnarray*}
X_t &= & \widetilde{B}_t+2\mu \int_0^t \left(\sum\limits_{j=0}^{i-1} \indi_{[\sigma_{2i},\sigma_{2i+1}[}(s) \frac{S(B_s^\mu)}{S(\epsilon)-S(B_s^\mu)}\right)ds + \mu t, \\
 &= & \widetilde{B}_{\sigma_{2i -1}}+2\mu \int_0^{\sigma_{2i -1}} \left(\sum\limits_{j=0}^{i-1} \indi_{[\sigma_{2i},\sigma_{2i+1}[}(s) \frac{S(B_s^\mu)}{S(\epsilon)-S(B_s^\mu)}\right)ds + \mu {\sigma_{2i -1}} \\
 &&+ (\widetilde{B}_t -\widetilde{B}_{\sigma_{2i -1}}) + \mu (t-{\sigma_{2i -1}}),\\
 &= & X_{\sigma_{2i -1}}+ (\widetilde{B}_t -\widetilde{B}_{\sigma_{2i -1}}) + \mu (t-{\sigma_{2i -1}}).
 \end{eqnarray*}
Since $X_{\sigma_{2i -1}} = B_{\sigma_{2i -1}}^\mu = -\alpha$, we have
$$ X_t = -\alpha + (\widetilde{B}_t -\widetilde{B}_{\sigma_{2i -1}}) + \mu (t-{\sigma_{2i -1}}).$$
Thus, $X$ is a Brownian motion with drift $\mu$ starting at time $\sigma_{2i -1}$.
  \item  If $t \in  [\sigma_{2i }, \sigma_{2i +1}[$,
  \begin{eqnarray*}
X_t &=& \widetilde{B}_t +  2\mu \int_0^t \left(\sum\limits_{j=0}^{i} \indi_{[\sigma_{2j},\sigma_{2j+1}[}(s) \frac{S(B_s^\mu)}{S(\epsilon)-S(B_s^\mu)}\right)ds + \mu t\\
&=& \widetilde{B}_{\sigma_{2i }} +  2\mu \int_0^{\sigma_{2i }}  \left(\sum\limits_{j=0}^{i-1} \indi_{[\sigma_{2j},\sigma_{2j+1}[}(s) \frac{S(B_s^\mu)}{S(\epsilon)-S(B_s^\mu)}\right)ds  + \mu {\sigma_{2i }} + \\
&& (\widetilde{B}_t-\widetilde{B}_{\sigma_{2i }}) + 2\mu \int_{\sigma_{2i }}^t  \frac{S(B_s^\mu)}{S(\epsilon)-S(B_s^\mu)}ds  + \mu (t-{\sigma_{2i }}),
 \\
&=&X_{\sigma_{2i }} + (\widetilde{B}_t-\widetilde{B}_{\sigma_{2i }}) +\mu \int^t_{\sigma_{2i }}
\left(\frac{2S(X_s)}{S(\epsilon)-S(X_s)}+1\right)ds.
\end{eqnarray*}
Since $X_{\sigma_{2i }} = 0$, we have
$$
X_t = (\widetilde{B}_t-\widetilde{B}_{\sigma_{2i }}) -\mu \int_{\sigma_{2i }}^t \coth(\mu(\epsilon-X_s))ds.$$
Heuristically, $X$ is the process $B^\mu$ conditioned not to hit the level $\epsilon$ (see \cite{a6}).
\end{itemize}
As a result, $X=B^\mu$ conditoned by $A_n$ is exactly the process $X$ defined by Definition \ref{definitionX}. Since $S= S_0 e^{\sigma B^\mu}$ and $Z= S_0 e^{\sigma X}$, the result holds.
\qed

\textbf{3. Complement.} We can construct a similar process $X$ by using the local time at 0. Instead of conditioning by the number of downcrossings, we condition by the amount of time spent at time 0.  The construction of $X$ and the computation of the strategy are almost similar to the ones in Section \ref{sec1}, \ref{sec1bis}, \ref{sec2} and \ref{sec3} and an explicit formula is given.

The asset $X$ has the following description. Let $\delta$ be a random variable with exponential law $\she(a)$. $\delta$ is the minimal time spent at 0 before being allowed to go through the resistance. Let $L_t$ be the local time at 0 of $X$ and  $ \tau:= \inf \{t \geqslant 0, L_t = \delta\}$.
\begin{itemize}
\item For all $t$ such that $X_t< 0$, or $t > \delta$, $X$ is a simple brownian motion with drift $\mu$.
\item For all $t$ such that $ X_t > 0$ and $L_t \leqslant \delta$ (i.e. we have not spent enough time at 0), we have:
$$
 X_t=\widetilde{B}_t+\mu \int_0^{t} \left(- \coth (\mu( \epsilon - X_s))  \indi_{ \{  X_s> 0 \}} +  \indi_{ \{ X_s \leqslant 0 \}}\right) ds,
$$
where $\widetilde{B}_t$ is a standard brownian motion.
\end{itemize}
Then, the optimal strategy for the logarithmic utility function is:
$$\pi_t^\star=\frac{\mu_0-r}{\sigma^2}
+\frac{2\mu  \frac{S(B_t^\mu)}{S(\epsilon) - S(0)}  \frac{a e^{-a L_t}}{\frac{\mu}{1-e^{-2\mu\epsilon}}  +a}  \indi_{ \{ B_t^\mu > 0 \}}  }{\sigma (1 -  e^{-aL_t})  +
\sigma  \frac{S(\epsilon) - S(B^\mu_t)}{S(\epsilon) - S(0)}  \frac{a e^{-a L_t}}{\frac{\mu}{1-e^{-2\mu\epsilon}}  +a} }.$$
The results look like to the one of Sections \ref{sec1}, \ref{sec2} and \ref{sec3}. However, such $X$ and such strategies can not be easily simulated because of the local time. This is why we have decided to focus on a model with downcrossings.

\section{Study in the case where the number of downcrossings is fixed}\label{sec2}

We consider $Z= S_0 e^{\sigma X}$ the price process defined in Section \ref{sec1bis}. In this Section, the minimal number of downcrossings is a fixed integer $n$ and we evaluate the optimal strategy for the logarithmic utility function. The main Theorem is in item \textbf{1} and its proof is in item \textbf{2}.

\medskip
\textbf{1.} We recall that $Z$ (resp. $X$) is $S$ (resp. $B^\mu$) conditioned to do at least $n$ downcrossings and that the expected wealth  for the logarithmic utility function is given by $ \E[\log(W^\pi_T)|A_n]$.
\begin{theo}
\label{stratgienfixe}
The strategy $\pi^\star$ which maximizes $\E[\log(W^\pi_T)|A_n]$ is
    $$ \pi_t^\star = \frac{ \mu_0 - r }{\sigma^2} -   \frac{2\mu}{\sigma}  \sum_{i=0}^{n-1} \indi_{t \in  [\sigma_{2i }, \sigma_{2i +1}[}
   \frac{S(X_t)}{S(\epsilon) - S(X_t)} ,
     $$
     for all $t \leqslant T_\epsilon$.
\end{theo}

\textbf{Comparison with the classic problem of section \ref{sec0}.} Without conditioning by the number of downcrossings, we have
$$\pi_s^\star  =  \frac{ \mu_0 - r }{\sigma^2}, \qquad \E[\log(W^{\pi^\star}_T)] = \log W_0 + rT + \frac{T (\mu_0 -r)^2}{2\sigma^2}.$$
When we condition by $A_n$, we get the result of Theorem \ref{stratgienfixe} and
$$
\E[\log(W^{\pi^\star}_T)| A_n] =\log W_0 + rT + \frac{T (\mu_0 -r)^2}{2\sigma^2}
$$
$$  + \frac{2 \mu}{p^n}   \int_0^T \E \left[ \sum_{i=1}^n \indi_{s \in  [\sigma_{2i }, \sigma_{2i +1}[} p^{n-1-i}
\frac{(\mu_0 -r) S(B_s^\mu)S(\epsilon)+ (\mu_0+r-\sigma^2)(S(B_s^\mu))^2}{\sigma(S(\epsilon) - S(B_s^\mu))(S(\epsilon) - S(-\alpha))}
 \right]ds
$$
We can see that, in comparison with the classic case, conditioning adds a new term  which only takes into account the downcrossings.

\medskip
\textbf{2. Proof.} The proof is divided into two steps. We start by computing explicitly $\E[\log(W^\pi_T)|A_n]$ in the first step. Then, we deduce the optimal strategy in the second one.

\textbf{a) Computing of $\E[\log(W^\pi_T)|A_n]$.}  We recall  that $W^\pi$ is solution of the equation (\ref{Wbis}).

We want to condition by $A_n$. Let us use Proposition \ref{Btilde} to replace $dB_t$ in (\ref{Wbis}):
\begin{equation}
\label{deterministe04}
\frac{d W^\pi_t}{W^\pi_t} = \alpha_0 (t) dt + \alpha_1(t) d\widetilde{B}_t,
\end{equation}
with
\begin{eqnarray*}
\alpha_0(t) &=& \pi_t \mu_0 + (1- \pi_t) r + 2 \mu \pi_t \sigma  \sum_{i=0}^{n-1} \indi_{t\in [\sigma_{2i}, \sigma_{2i + 1} [} \frac{S(B_t^\mu)}{S(\epsilon) - S(B_t^\mu)}, \\
\alpha_1(t) &=& \pi_t \sigma.
\end{eqnarray*}
Via It\^o's formula, we deduce that the solution of the equation (\ref{deterministe04}) is:
$$ W_t^\pi = W_0^\pi \exp \left( \int_0^t \alpha_1(s) d\widetilde{B}_s + \int_0^t \left(\alpha_0(s)- \frac{\alpha_1^2(s)}{2} \right) ds \right).$$
Since $\widetilde{B}$ is independent from $A_n$, we get
$$
 \E[\log(W^\pi_T) \indi_{A_n}] =
 \log ( W_0^\pi) \Pb(A_n) + \E\left[ \int_0^T \alpha_1(s) d\widetilde{B}_s  \right] \Pb(A_n)
 + \E\left[ \int_0^T \left(\alpha_0(s)- \frac{\alpha_1^2(s)}{2} \right)  \indi_{A_n} ds \right]
 $$
The expressions of $\alpha_0(s)$ and $\alpha_1(s)$ give then
\begin{eqnarray*}
 \E[\log(W^\pi_T) \indi_{A_n}]  &= &( \log ( W_0^\pi) + r T )\Pb(A_n) \\
 &&+  \E\left[   \indi_{A_n}  \int_0^T \pi_s  \left(\mu_0 - r + 2 \mu \sigma  \sum_{i=0}^{n-1} \indi_{s\in [\sigma_{2i}, \sigma_{2i + 1}[} \frac{S(B_s^\mu)}{S(\epsilon) - S(B_s^\mu)} \right)  - \frac{\pi_s^2 \sigma^2}{2} ds \right] .
\end{eqnarray*}
Finally, we obtain:
\\

$\displaystyle \E[\log(W^\pi_T)|A_n]  =  \log ( W_0^\pi)+ r T $
\begin{equation}\label{Wdesfixecalcul}+ \frac{1}{\Pb(A_n)} \E\left[   \indi_{A_n}  \int_0^T \pi_s  \left(\mu_0 - r + 2 \mu \sigma  \sum_{i=0}^{n-1} \indi_{s\in [\sigma_{2i}, \sigma_{2i + 1}[} \frac{S(B_s^\mu)}{S(\epsilon) - S(B_s^\mu)}   - \frac{\pi_s \sigma^2}{2} \right)  ds \right].
\end{equation}

\textbf{b) Computation of the optimal strategy $\pi^\star $.} To maximize $\E[\log(W^\pi_T)|A_n] $, we have to maximize the term in the integral, for all $s\in [0,T]$ :
\begin{equation}
\label{deterministe05}
 \pi_s  \left(\mu_0 - r + 2 \mu \sigma  \sum_{i=0}^{n-1} \indi_{s\in [\sigma_{2i}, \sigma_{2i + 1}[} \frac{S(B_s^\mu)}{S(\epsilon) - S(B_s^\mu)} \right)    - \frac{\pi_s^2 \sigma^2}{2} .
\end{equation}
Let us distinguish 2 cases:
\begin{itemize}
  \item If $s \in  [\sigma_{2i -1}, \sigma_{2i }[$ (case where we are not in a downcrossing), then (\ref{deterministe05}) simplifies to
  $$  \pi_s  \left(\mu_0 - r  \right)    - \frac{\pi_s^2 \sigma^2}{2} \texttt{},$$
  and the maximum is reached for $\pi^\star_s = \frac{\mu_0 - r}{\sigma^2}$. Let us note that we recognize the classic formula.
  \item
  If $s \in  [\sigma_{2i }, \sigma_{2i +1}[$ (case where we are in a downcrossing), then (\ref{deterministe05}) becomes
  $$ \pi_s  \left(\mu_0 - r + 2 \mu \sigma \frac{S(B_s^\mu)}{S(\epsilon) - S(B_s^\mu)} \right)    - \frac{\pi_s^2 \sigma^2}{2} .$$
  The maximum is then reached for
  $$ \pi_s^\star = \frac{1}{\sigma^2} \left( \mu_0 - r + 2 \mu \sigma  \frac{S(B_s^\mu)}{S(\epsilon) - S(B_s^\mu)}  \right).$$
\end{itemize}
\qed


\section{Study in the case where the number of downcrossings is random}\label{sec3}

We consider $Z= S_0 e^{\sigma X}$ the price process defined in Section \ref{sec1bis}. In this Section, the minimal number of downcrossings is $\xi$  a random variable on $\N$, independent from $B^\mu$. We use the notation \(\alpha_n = \mathbb{P}(\xi = n)\). We compute an optimal strategy in item  \textbf{1}, without hypothesis on $\xi$ (The proof is postponed in item  \textbf{5}).  Then we consider particular laws for $\xi$: a geometric random variable in item  \textbf{2}, a random variable with compact support in item \textbf{3} and a variable with geometric queue in \textbf{4}.

 \textbf{1.}   In this section, $Z$ (resp. $X$) is $S$ (resp. $B^\mu$) conditioned to do at least $\xi$ downcrossings. We consider the event:
$$A_\xi:=
\begin{cases}
\Omega & \text{ for }\; \xi=0\\
\{\sigma_{2n-1}< T_\epsilon\} &\text{ for }\; \xi=n\geq1
\end{cases}
$$
i.e.  ``The trajectories $\omega$ which have done at least $\xi(\omega)$ downcrossings before hitting the level $\epsilon$". Then, the event $A_\xi$  has probability:
 $$\Pb(A_\xi)=\sum\limits_{n\geq0}\Pb(A_\xi\cap\{\xi=n\})=\sum\limits_{n\geq0}\alpha_n\Pb(A_n)=\sum\limits_{n\geq0}\alpha_n p^n.$$
Then, the expected wealth  for the logarithmic utility function is given by $ \Q[\log(W^\pi_T)]$, where
$$\Q(\cdot):=\Pb[\cdot|A_\xi].$$
The optimal strategy for $Z$ is given by the following theorem.
\begin{theo}
\label{theostrategiedescent}
The strategy $\pi^\star$ which maximizes $\Q[\log(W^\pi_T)]$ is :
$$\pi_t^\star=\frac{\mu_0-r}{\sigma^2}+\frac{2\mu  \frac{S(B_t^\mu)}{S(\epsilon)-S(B_t^\mu)} \sum\limits_{n\geq0} \alpha_n M_t^n\sum\limits_{i=0}^{n-1}\indi_{[\sigma_{2i},\sigma_{2i+1}[}(t)}{\sigma \sum\limits_{n\geq0} \alpha_n M_t^n},$$
for $t < T_\epsilon$.
\end{theo}

Following the optimal strategy, we obtain for maximum:
\begin{equation}
\label{maxdescent}
\Q[\log(W_T^{\pi^\star})]=\log(W_0)+rT+\frac{1}{\Pb(A_\xi)}\int_0^T\E\left[\frac{\sigma^2}{2}\sum\limits_{n\geq0} \alpha_n M_t^n \left(\pi_t^\star\right)^2  \right]dt.
\end{equation}

\begin{remark}
$\pi^\star$ is composed of two terms: the first one, $\frac{\mu_0-r}{\sigma^2}$, corresponds to the optimal allocation strategy for the classic problem, when the underlying process is a geometric Brownian motion. The second term appears only during downcrossings. It is always negative, which suggests that we should reduce our proportion of risky assets in a downcrossing.
\end{remark}

\textbf{In an upcrossing:} for $t\in [\sigma_{2i_0-1}, \sigma_{2i_0}[$ and $t<T_\epsilon$, the expression of $\pi^\star$ is the classic formula:
$$\pi_t^\star=\frac{\mu_0-r}{\sigma^2}.$$
And we have
\begin{equation}\label{maxdescentmgen}
\sum\limits_{n=0}^\infty\alpha_n M_t^n
=\sum\limits_{n=0}^{i_0}\alpha_n + \sum\limits_{n=i_0+1}^\infty \alpha_n p^{n-i_0}.
\end{equation}

\textbf{In a downcrossing:} for $t\in [\sigma_{2i_0}, \sigma_{2i_0+1}[$ and $t<T_\epsilon$, we have
\begin{equation}
\label{maxdescentdgen}
\sum\limits_{n=0}^\infty\alpha_nM_t^n
=\sum\limits_{n=0}^{i_0}\alpha_n + \sum\limits_{n=i_0+1}^\infty \alpha_n p^{n-1-i_0} \frac{S(\epsilon)-S(B_t^\mu)}{S(\epsilon)-S(-\alpha)},
\end{equation}
and:
\begin{equation}
\label{strategiedescentgen}
\pi_t^\star=\frac{\mu_0-r}{\sigma^2}+\frac{  2\mu   \frac{S(B_t^\mu)}{S(\epsilon)-S(-\alpha)} \sum\limits_{n=i_0+1}^\infty  \alpha_n p^{n-1-i_0} } {\sigma\left(\sum\limits_{n=0}^{i_0}\alpha_n + \frac{S(\epsilon)-S(B_t^\mu)}{S(\epsilon)-S(-\alpha)} \sum\limits_{n=i_0+1}^\infty \alpha_n p^{n-1-i_0} \right)}.
\end{equation}

\textbf{2. Case of $\xi$ being a geometric random variable.} In this item, the law of $\xi$ is the following:
$$\Pb(\xi=n):=\alpha_n=q^n(1-q), \textrm{ where } q \in [0,1[.$$
Then, the event $A_\xi$  has probability:
$$\Pb(A_\xi)=\sum\limits_{n\geq0}q^n(1-q)p^n=\frac{1-q}{1-pq}.$$
During downcrossings and upcrossings, the expression of $\pi^\star$ and $\sum\limits_{n=0}^\infty\alpha_nM_t^n$  can be simplified.

\textbf{In an upcrossing:} for $t\in [\sigma_{2i_0-1}, \sigma_{2i_0}[$ and $t<T_\epsilon$, the expression of $\pi^\star$ is the classic formula: $\pi_t^\star=\frac{\mu_0-r}{\sigma^2}$ and we have
$$\sum\limits_{n=0}^\infty\alpha_nM_t^n = 1-q^{i_0+1} + \frac{(1-q)pq^{i_0+1}}{1-pq}.
$$

\textbf{In a downcrossing:} for $t\in [\sigma_{2i_0}, \sigma_{2i_0+1}[$ and $t<T_\epsilon$, we have
$$
\sum\limits_{n=0}^\infty\alpha_nM_t^n
= 1-q^{i_0+1} + \frac{S(\epsilon)-S(B_t^\mu)}{S(\epsilon)-S(-\alpha)} \frac{(1-q)q^{i_0+1}}{1-pq}.
$$
As a result, we get:
\begin{equation}
\label{strategiedescent}
\pi_t^\star=\frac{\mu_0-r}{\sigma^2}+\frac{2\mu\frac{S(B_t^\mu)}{S(\epsilon)-S(-\alpha)} \frac{(1-q)q^{i_0+1}}{1-pq}}{\sigma\left(1-q^{i_0+1} + \frac{S(\epsilon)-S(B_t^\mu)}{S(\epsilon)-S(-\alpha)} \frac{(1-q)q^{i_0+1}}{1-pq}\right)}.
\end{equation}

\textbf{3. Case of $\xi$ having compact support.} Instead of allowing $\xi$ to have any value, we consider that $\{ \xi=n \}$ for huge values of $n$ has no meaning. As a matter of fact, a trajectory which has a huge number of downcrossings before getting to the resistance level is not likely to appear in price charts. Thus, we consider that $\alpha_n= 0$ for all $n > N$, with $N$ fixed in $\N$. As a result, $\Pb(A_\xi)= \sum\limits_{n=0}^N \alpha_n p^n$ and  the sums in (\ref{maxdescentdgen}), (\ref{maxdescentmgen}) and (\ref{strategiedescentgen}) are finite sums.

Let us consider that $i_0 > N$ and $t<T_\epsilon$, i.e. we have more than $N$ downcrossings.  Then, for $t \geqslant \sigma_{2i_0-1}$, we have
$$
 \sum\limits_{n=0}^\infty\alpha_n M_t^n =\sum\limits_{n=0}^{N}\alpha_n = 1 $$
and
$$
\pi_t^\star=\frac{\mu_0-r}{\sigma^2}.
$$
We are back to the classic strategy.

\textbf{4. Case of $\xi$ having a geometric queue.} Although the case $\xi$ geometric (item \textbf{2}) gives nice mathematical formulas, this case weights too much the event $\xi=0$, that is the classic case without forced downcrossings. On the other hand, the case of item \textbf{3} forgets very indecisive patterns, where the trajectory has lots of downcrossings. By mixing the two laws, we can have a law which weights more on the $N^{th}$ first downcrossings without erasing huge numbers of downcrossings.

Let $\alpha_0, \dots, \alpha_N$ be such that $\sum_{i=0}^N \alpha_n = S \in ]0,1[$. We put 
 \(q = (1-S)^{1/(N+1)}\).
for $n> N$. Then, $\sum_{n=0}^\infty \alpha_n = 1$ and $\xi$ is a random variable with geometric queue.

Then, $A_\xi$ has probability:
 $$\Pb(A_\xi)= \sum\limits_{n=0}^N \alpha_n p^n +  \sum\limits_{n\geq N+1}q^n(1-q) p^n= \Pb(A_\xi)= \sum\limits_{n=0}^N \alpha_n p^n + (pq)^{N+1} \frac{1-q}{1-pq}.$$

 During downcrossings and upcrossings, the infinite sums in expression of $\pi^\star$ and $\sum\limits_{n=0}^\infty\alpha_nM_t^n$  can be calculated.

\textbf{In an upcrossing:} for $t\in [\sigma_{2i_0-1}, \sigma_{2i_0}[$ and $t<T_\epsilon$, we have
$$\sum\limits_{n=0}^\infty\alpha_nM_t^n = \sum\limits_{n=0}^{i_0}\alpha_n  + \sum\limits_{n=i_0+1}^{N}\alpha_n p^{n-i_0}  + \frac{(1-q)p^{(N+1-i_0)\vee 1}q^{(N\vee i_0)+1}}{1-pq}.
$$
with the convention $\sum_a^b=0$ if $b<a$.

\textbf{In a downcrossing:} for $t\in [\sigma_{2i_0}, \sigma_{2i_0+1}[$ and $t<T_\epsilon$, we have
$$
\sum\limits_{n=0}^\infty\alpha_nM_t^n
= \sum\limits_{n=0}^{i_0}\alpha_n  + \frac{S(\epsilon)-S(B_t^\mu)}{S(\epsilon)-S(-\alpha)} \left(  \sum\limits_{n=i_0+1}^{N}\alpha_n p^{n-1-i_0} +  \frac{(1-q)q^{(N\vee i_0)+1} p^{(N-i_0)^+}}{1-pq}\right).
$$
As a result, we get:
$$
\pi_t^\star=\frac{\mu_0-r}{\sigma^2}+\frac{  2\mu   \frac{S(B_t^\mu)}{S(\epsilon)-S(-\alpha)}
\left(  \sum\limits_{n=i_0+1}^N\alpha_n p^{n-1-i_0} + q^{(N\vee i_0)+1}p^{(N-i_0)^+} \frac{1-q}{1-pq} \right) }{\sigma\left(\sum\limits_{n=0}^{i_0}\alpha_n  + \frac{S(\epsilon)-S(B_t^\mu)}{S(\epsilon)-S(-\alpha)} \left(  \sum\limits_{n=i_0+1}^{N}\alpha_n p^{n-1-i_0} +  \frac{(1-q)q^{(N\vee i_0)+1} p^{(N-i_0)^+}}{1-pq}\right) \right)}.
$$

\textbf{5. Proof.} We start by computing  $\Q[\log(W_T^\pi)]$. The key of this step is to introduce $A_n$. We make the decomposition:
$$ \E\left[\log(W_T^\pi) \indi_{A_\xi}\right]= \sum\limits_{n\geq0} \E\left[\log(W_T^\pi) \indi_{A_n}|\xi=n\right]\alpha_n =\sum\limits_{n\geq0} \E\left[\log(W_T^\pi) \indi_{A_n}\right]\alpha_n.$$
Using (\ref{Wdesfixecalcul}) of Section \ref{sec2}, we get:
\begin{multline*}
\E\left[\log(W_T^\pi)\indi_{A_\xi}\right]=  (\log(W_0)+rT) \sum_{n \geqslant 0}\alpha_n\Pb(A_n)\\
+ \sum\limits_{n\geq0} \alpha_n \E\left[\indi_{A_n}\int_0^T\pi_t\left(\mu_0-r+2\mu\sigma\sum\limits_{i=0}^{n-1}\indi_{[\sigma_{2i},\sigma_{2i+1}[}(t)\frac{S(B_t^\mu)}{S(\epsilon)-S(B_t^\mu)}-\frac{\sigma^2}{2}\pi_t\right)dt\right].
\end{multline*}
In the expectation, we condition by $\shf_t$. Since $\Pb(A_n| \shf_t)= M_t^n$, it comes:
\begin{multline*}
\E\left[\log(W_T^\pi)\indi_{A_\xi}\right]=  \alpha_n (\log(W_0)+rT) \Pb(A_\xi)\\
+ \sum\limits_{n\geq0} \alpha_n \int_0^T \E\left[M_t^n\pi_t\left(\mu_0-r+2\mu\sigma\sum\limits_{i=0}^{n-1}\indi_{[\sigma_{2i},\sigma_{2i+1}[}(t)\frac{S(B_t^\mu)}{S(\epsilon)-S(B_t^\mu)}-\frac{\sigma^2}{2}\pi_t\right)\right] dt.
\end{multline*}
Therefore, we finally get
\begin{multline*}
 \Q[\log(W_T^\pi)] = \log(W_0)+rT \\
+ \frac{1}{\Pb(A_\xi)}\int_0^T \E\left[\pi_t \sum\limits_{n\geq0} \alpha_n M_t^n\left(\mu_0-r+2\mu\sigma\sum\limits_{i=0}^{n-1}\indi_{[\sigma_{2i},\sigma_{2i+1}[}(t)\frac{S(B_t^\mu)}{S(\epsilon)-S(B_t^\mu)}-\frac{\sigma^2}{2}\pi_t\right)\right] dt.
\end{multline*}
The optimal allocation strategy is found by maximizing the term in brackets:
\[
\pi_t^\star=\frac{\mu_0-r}{\sigma^2}+\frac{2\mu  \frac{S(B_t^\mu)}{S(\epsilon)-S(B_t^\mu)} \sum\limits_{n\geq0} \alpha_n M_t^n\sum\limits_{i=0}^{n-1}\indi_{[\sigma_{2i},\sigma_{2i+1}[}(t)}{\sigma \sum\limits_{n\geq0} \alpha_n M_t^n}.
\]
\qed

\section{Simulation and comparison with the classic results}\label{sec4}

The aim of this numerical part is to compare our strategy with the classical one. We simulate a trajectory of $X$ which does  a random number of downcrossings and we compute the optimal strategy and the resulting wealth.  In order to simulate this kind of trajectory, we must take into account the two distinct phases of the movement.
\begin{itemize}
\item During an upcrossing, $X$ is a standard Brownian motion with drift $\mu$. $\pi^\star$ is the classic strategy (\ref{strategieinitial}).
\item During a downcrossing, $X$ follows the SDE (\ref{xdescent}). $\pi^\star$ is calculated using (\ref{strategiedescent}).
\end{itemize}
%
%
%
%
%
%
%
%
%
%
Since $\pi$ is the proportion of the risky asset, we impose $\pi$ to stay in $[0,1]$. Then, the optimal strategy for the constrained optimization problem is the optimal strategy of the unconstrained problem, projected on $[0,1]$. Thus, we project  $\pi$ and $\pi^\star$ on $[0,1]$. At the same time, we evaluate the wealth created by following the classic strategy $\pi^c$ (\ref{strategieinitial}), in order to compare the efficiency of $\pi^\star$ with the one of $\pi^c$.

We fix the following parameters: $ \mu_0 =0,1$, $  \sigma = 0,15$, $  r = 0,02 $, $ \alpha =1 $ and $  \epsilon = 0,3 $. Then, $\mu = 0,5916 $ and $p=0,1165$. With our parameter, $ \pi^c$ projected on $[0,1]$ is equal to 1 all the time.

\textbf{1.} First, we consider that $\xi$ is a geometric random variable (c.f. item \textbf{2} of section \ref{sec3}) with parameter $q$. Under $A_\xi$, the probability to have no forced downcrossing is
$$\Pb(\xi = 0 | A_\xi)= \frac{\Pb( \xi=0 \cap A_\xi)}{\Pb{A_\xi}} =  \frac{\Pb( \xi=0) \Pb ( A_0)}{\Pb{A_\xi}}= 1-pq.$$
Since $p= 0,1165$, the probability to have no forced downcrossing is at least $0,88$.

No forced downcrosssing means that $X$ is the Brownian motion with drift all the time and there is no resistance level. Note that our model consider that we are in a downcrossing since the beginning and thus, we do not follow the optimal strategy (which is $\pi^c$). The case without forced downcrossing is an unfavorable case for us. Trajectories without resistance are not the topic of our study. Since the probability of this event is too high, we consider that choosing $\xi$ geometric is not relevant for simulation.

\textbf{2.} Then, we consider that  $\xi$ has finite support $\{0, \dots, N\}$.

On one side, knowing the law of $\xi$ under $\Pb$ (i.e $\alpha_n, n \geq 0$), it is easy to deduce the law of $\xi$ under $\Q$:
$$ \beta_n = \Q( \xi=n)= \Pb( \xi = n | A_\xi) = \frac{\Pb( \xi = n \cap A_\xi)}{\Pb(A_\xi)} =
\frac{\Pb( \xi = n) \Pb( A_n)}{\Pb(A_\xi)} =\frac{\alpha_n p^n}{\sum_{i=0}^N \alpha_i p^i}.$$
On the other side, the observations on real charts are observations of $\xi$ under $\Q$, i.e we can calibrate only the values of $\beta_n$. However, we need the values of $\alpha_n$ in order to compute $\pi^*$. Then, it is necessary to be able to compute $\alpha_n, n \in [0,N],$ thanks to $\beta_n, n \in [0,N]$. In the case of finite support, we have an explicit formula:
$$ \alpha_n = \frac{\beta_n}{p^n \sum_{i=0}^N \frac{\beta_i}{p^i}}.$$
A resistance line can be drawn only if there are at least three aligned local maxima. However $\beta_0$ and $\beta_1$ are not necessary 0, since our model starts after these three maxima and the drawing of the resistance line.

In our simulation, we fix $N=6$ and
\[
 \begin{array}{rclrclrcl}
\beta_0 &=& 0,1  ,& \beta_1 &=& 0,1   ,& \beta_2 &=& 0,2  ,\\
\beta_3 &=& 0,2  ,&  \beta_4 &=& 0,2  ,& \beta_5 &=& 0,1  ,\\
\beta_6 &=& 0,1  .& &&   & &&
\end{array}
\]
We simulate $X, \pi^\star, \pi^c, W^{\pi^\star},W^{\pi^c}$ with an initial wealth $W_0=1$.  Let us present a typical result of the simulation in the following figures.

\begin{tabular}{c c}
\begin{tabular}{p{7cm}}
\textit{ See the jointed picture.
 Latex Compilation via HAL
  does not allow jpg. 
  Thus, we had to remove
   it from our article.}
\end{tabular}
&
\begin{tabular}{p{6cm}}
The first figure shows a trajectory of $X$ which does three downcrossings and the strategy $\pi^\star$ at the same time. \\
\hfill\\
The second figure compares the wealth $W^{\pi^\star}$ earned by following the strategy $\pi^\star$ with the wealth $W^{\pi^c}$ earned by following the classical strategy $\pi^c$. \\
\hfill\\
As the figure shows, the gain is better with the strategy $\pi^\star$, thanks to the downcrossings.
\end{tabular}
\end{tabular}
During downcrossings, $\pi^\star$ (projected on $[0,1]$) vanishes, which means that we invest all our wealth in the bond. Since we expect the risky asset to fall down, the more natural way to preserve our wealth is to secure it in the bond. Thus, the wealth rises slowly during ``bad days" (downcrossing).  On the contrary, the standard strategy keeps investing in the risky asset and $W^{\pi^c}$ decreases. Consequently, at the end of a downcrossing,  $W^{\pi^\star}$ is likely to be higher than $W^{\pi^c}$.

When the level $-\alpha$ is reached, the forced downcrossing ends and the price can go freely. Then, we expect the risky asset to rise. Following the strategy  $\pi^\star$, we invest in the risky asset and the wealth will likely follows the rise of the risky asset. The strategy $W^{\pi^c}$ does the same, but since the wealth $W^{\pi^\star}$ at the beginning of the rise is higher than $W^{\pi^c}$, we still gain more during the rise that the standard strategy.

When the resistance level 0 is reached, we expect to have a fall of the price and we secure our wealth in the bond. This process last until we have done three downcrossings. After this time, $X$ is not anymore constrained between the levels 0 and $\alpha$. Thus, we are back to a classical Brownian motion and a classical strategy.

\medskip\textbf{3.} The second part of our simulation is to compare the wealths at the final time $T=10$. Via Monte-carlo methods, we evaluate $W^{\pi^c}_T-W^{\pi^\star}_T$. The mean value is $0,8415$ with a standard deviation of $0,005$ for $2000$ simulations.

The two following tables shows the results when we change only the volatility $\mu_0$ (first table), and only the level $\alpha$ (second table).
\[
 \begin{array}{|c|c|c|}
 \hline
 \mu & \textrm{mean} & \textrm{std dev.}  \\
 \hline
 \hline
0,1 &  0,8415  & 0,005 \\
  \hline
 0,15 & 1,1347 & 0,0095 \\
  \hline
 0,2 & 1,6394 & 0,0185\\
 \hline
\end{array}
\]
It seems that the more volatile the price is, the more we earn. But the precision of our evaluation is lower.

The second table shows the influence of $\alpha$. We have done only 1000 simulations for theses results.
\[
 \begin{array}{|c|c|c|}
 \hline
 \alpha & \textrm{mean} & \textrm{std dev.}  \\
 \hline
 \hline
0.5&  0.4889  & 0.0051 \\
  \hline
 0.6 & 0.5581 & 0.0066 \\
  \hline
 0.7 & 0.6238 & 0.0078\\
  \hline
0.8 & 0.6846 & 0.0086 \\
 \hline
0.9 & 0.7477 & 0.01 \\
 \hline
1 & 0.8618  & 0.0112 \\
 \hline
1.1 & 0.927 & 0.0137\\
 \hline
1.2 & 1.0139 &0.0177 \\
 \hline
\end{array}
\]
The mean gain grows with $\alpha$. $\alpha$ represents the ``minimum height of downcrossing" after reflecting on the resistance. Higher values of $\alpha$ implies that the fall of price is more important. Consequently, the loss with the standard strategy is higher.

\section{Conclusion}
In the article, we have introduced a variation on standard Black-Scholes model.
Our main purpose was to propose a stochastic model with paths having properties
similar to observed prices. Our model have resistance (or support) lines.
We applied well-known results of stochastic analysis to this new model and
exhibited an optimal allocation portfolio procedure.

Many other rules issued from technical analysis in finance are mysterious 
from a financial mathematics point of view. It should be interesting
to improve the accuracy of the mathematical models to the observations.

\section*{Acknowledgement}
The authors would like to thank gratefully Pr. Pierre Vallois for fruitful discussions.

\bibliography{bibfinance}
\bibliographystyle{plain}
\nocite{*}

\end{document}